\documentclass[conference]{IEEEtran}
\IEEEoverridecommandlockouts
\usepackage{cite}
\usepackage{amsmath,amssymb,amsfonts}
\usepackage{amsthm}
\usepackage{algorithmic}
\usepackage{graphicx}
\usepackage{textcomp}
\usepackage{xcolor}

\usepackage{multirow}
\usepackage{booktabs}
\usepackage[algo2e,ruled,linesnumbered]{algorithm2e}

\def\BibTeX{{\rm B\kern-.05em{\sc i\kern-.025em b}\kern-.08em
    T\kern-.1667em\lower.7ex\hbox{E}\kern-.125emX}}
\begin{document}

\title{Conformal Uncertainty Quantification of Electricity Price Predictions for Risk-Averse Storage Arbitrage  \\

\thanks{The authors would like to acknowledge supports from Columbia Data Science Institute, Red River Clean Energy, King Abdulaziz City for Science and Technology (KACST), and National Science Foundation under award ECCS-2239046.}
}

\author{\IEEEauthorblockN{Saud Alghumayjan, Bolun Xu}
\IEEEauthorblockA{\textit{Dept. of Earth and Environmental Engineering} \\
\textit{Columbia University}\\
New York, NY 10027, USA \\
\{saa2244, bx2177\}@columbia.edu}
\and
\IEEEauthorblockN{Ming Yi}
\IEEEauthorblockA{\textit{Data Science Institute} \\
\textit{Columbia University}\\
New York, NY 10027, USA \\
my2826@columbia.edu}
}

\maketitle

\begin{abstract}
This paper proposes a risk-averse approach to energy storage price arbitrage, leveraging conformal uncertainty quantification for electricity price predictions. The method addresses the significant challenges posed by the inherent volatility and uncertainty of real-time electricity prices, which create substantial risks of financial losses for energy storage participants relying on future price forecasts to plan their operations. The framework comprises a two-layer prediction model to quantify real-time price uncertainty confidence intervals with high coverage. The framework is distribution-free and can work with any underlying point prediction model. We evaluate the quantification effectiveness through storage price arbitrage application by managing the risk of participating in the real-time market. We design a risk-averse policy for profit-maximization of energy storage arbitrage to find the safest storage schedule with very minimal losses. Using historical data from New York State and synthetic price predictions, our evaluations demonstrate that this framework can achieve good profit margins with less than $35\%$ purchases.
\end{abstract}

\begin{IEEEkeywords}
Electricity markets, Energy storage, Machine learning
\end{IEEEkeywords}

\section{Introduction}\label{sec:introduction}

Rapid integration of renewable resources introduced unprecedented levels of operational uncertainty and reliability challenges, especially in real-time power system operations. These factors materialized into volatile electricity prices~\cite{wang2017impact}, which attract surging energy storage participation in wholesale markets~\cite{us2023data}, charging their systems during low-price periods and discharging during price peaks. The real-time market holds particular appeal for storage participants, as it offers greater price volatility than other market options and, consequently, higher potential profits~\cite{qin2023role}.

However, optimal operation requires accurate forecasting of future prices, which has become increasingly challenging due to multiple factors, including real-time operational uncertainty and demand forecasting errors. Researchers have applied a variety of models to enhance real-time price forecasting, utilizing techniques that range from statistical methods to advanced artificial intelligence (AI) frameworks~\cite{lago2018forecasting}. Most of the research on electricity price prediction focuses on point forecasting, with models like Long Short-Term Memory (LSTM), Multi-Layer Perceptron (MLP)~\cite{panapakidis2016day}, and hybrid models~\cite{lehna2022forecasting} being popular due to their strong capability to capture complex patterns in historical prices. Moreover, models from other areas of machine learning have also been adapted to time-series problems, such as electricity price and transformer variant models~\cite{alghumayjan2024energy}. More recently, ideas for decision-focused frameworks such as the predict-then-optimize approach, which incorporates operational costs into the prediction methodology through hybrid loss functions~\cite{elmachtoub2022smart}. While this integration of prediction and optimization provides valuable insights, this approach does not directly address uncertainty quantification in price predictions, which remains an important consideration for managing operational risks.

One of the prominent techniques for uncertainty quantification is the conformal prediction~\cite{vovk1999machine}. Where the idea is to construct prediction sets with certainly guaranteed coverage through a nonconformity score that measures how far the prediction is from the ground truth. Such guarantee is only achieved if the data is exchangeable~\cite{vovk2005algorithmic}. However, in time-series settings, the distribution shifts over time, while a conformal control prediction model for time series can achieve long-run coverage~\cite{angelopoulos2024conformal}. Conformal prediction has been explored in recent studies in day-ahead electricity prices~\cite{renkema2024conformal,yeh2024end}, while its application in the more volatile real-time markets, which provides higher profit expectations but also higher risk, has not been well-studied.

Real-time price arbitrage presents a more complex challenge due to the sequential nature of market clearing~\cite{zheng2022arbitraging}, making the integration of risk-averse objectives less straightforward~\cite{sioshansi2021energy}. This study bridges this gap by introducing a risk-averse uncertainty quantification framework based on distribution-free conformal prediction, offering compatibility with any underlying point prediction model.
Our key contributions are as follows:
\begin{itemize}
    \item We design a risk-averse policy for real-time energy storage arbitrage leveraging conformal uncertainty quantification for electricity price predictions.
    \item We build a conformal prediction model for multi-hour ahead real-time price uncertainty quantification with long-run coverage.
    \item We validate the proposed framework on historical data from New York State with synthetic price predictions. Our results show that our framework achieves very low losses while maintaining a good margin of profits.
\end{itemize}

The remainder of the paper is organized as follows: Section II introduces the formulation and methodology. Section III case study dataset, and Section IV concludes the paper.

\section{Formulation and Methodology}\label{sec:methodology}

\begin{figure*}[htbp]
\centerline{\includegraphics[trim = 0mm 0mm 0mm 0mm, clip, width=0.9\linewidth]{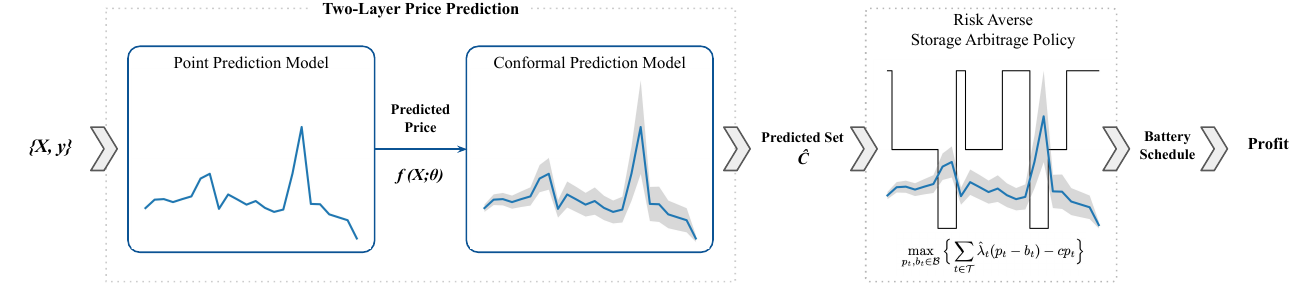}}
\caption{The pipeline for the proposed risk-averse storage arbitrage framework.}
\label{fig:conformal_diagram}
\end{figure*}

\subsection{Electricity Price Conformal Prediction}\label{sec:conformal}
We present a framework to predict and quantify the uncertainty in electricity real-time prices. The framework consists of two layers: the prediction layer and the conformal layer.

\textbf{The prediction layer}. We construct a point forecasting model employing machine learning techniques:
\begin{equation}
    \hat{\mathbf{\lambda}}_{t+1:t+T} = f(\mathbf{X}; \theta) \label{eq:point_model}
\end{equation}

In this model, $f(\mathbf{X}; \theta)$ represents any point forecast model, where $\theta$ denotes the model parameters, and $\mathbf{X}$ is the input multivariate time series. The input is defined as $ \mathbf{X} = [\mathbf{x}_{t-L+1}, \mathbf{x}_{t-L+2}, \ldots, \mathbf{x}_{t}]$, with $ \mathbf{x}_t \in \mathbb{R} $, and $L$ being the length of the look-back window. The model outputs a univariate time series forecast $\hat{\lambda}_{t+1}, \hat{\lambda}_{t+2}, \ldots, \hat{\lambda}_{t+T}$, where $T$ is the forecasting horizon. For calibration and testing, we divide the dataset into two subsets: $D_1$ for calibrating the point forecast predictions through the conformal control layer and $D_2$ for testing the calibrated model.
\begin{subequations}
\begin{gather}
     D_1 = \{(x_i,\lambda_i)\}_{i=1}^{T_1},\\
     D_2 = \{(x_i,\lambda_i)\}_{i=T_1+1}^{T_1+T_2}
\end{gather}
\end{subequations}
Where $T_1, \text{and } T_2$ represent the cutoffs for calibration and testing, respectively.
\\

\textbf{The conformal layer}. The goal here is to quantify the uncertainty of the price point prediction model by producing prediction sets that reflect how much the model is confident about its predictions. We follow a conformal prediction technique by using a portion of the data, the calibration set, to produce a $1-\alpha$ coverage that contains the ground truth. We construct the prediction sets through a nonconformity score $S(\lambda, \hat{\lambda})$, such as prediction signed residual \eqref{eq:score} that measures how far the prediction is from ground truth.
\begin{equation}
    S(\lambda, \hat{\lambda}) = \lambda - \hat{\lambda} \label{eq:score}
\end{equation}

The conformal confidence set is mathematically defined as~\eqref{eq:prediction_set}, where $q$ represents an estimate of the $1-\alpha$ quantile for the nonconformity score.

\begin{equation}
    \hat{C} = \{\lambda \in R: s(\lambda,\hat{\lambda}) \leq q\} \label{eq:prediction_set}
\end{equation}

As stated earlier in section~\ref{sec:introduction}, conformal prediction does not guarantee coverage if the data is not exchangeable. And electricity prices' distribution shifts over time. So, our goal here is to achieve long-run coverage ($\frac{1}{T} \sum_{t=1}^{T} 1\{\hat{\lambda}\notin \hat{C}\} = \alpha + o(1)$). The main idea is to construct the prediction sets by estimating the next quantile based on the PID feedback idea~\cite{angelopoulos2024conformal}:
\begin{equation}
    q_{t+1} = \hat{q}_{t+1} + r_t \sum_{i=1}^t (1\{\hat{\lambda}\notin \hat{C}\} - \alpha) \label{eq:control}
\end{equation}

where $\hat{q}_{t+1}$ is the quantile forecast to capture possible trends in the scores. $r_t$ is a saturation function that must follow $x \geq c \cdot h(t) \Rightarrow r_t(x) \geq b, \text{ and } x \leq -c \cdot h(t) \Rightarrow r_t(x)\leq-b$, for constants $b,c\geq0$, and a sub-linear, non-negative, non-decreasing function $h$. Angelopoulos et al. proved that \eqref{eq:control} achieve the long-run coverage, with any forecaster for $\hat{q}_{t+1}$ and any possible $r_t$, according to Theorem 1 \cite{angelopoulos2024conformal}.

\subsection{Real-Time Energy Storage Arbitrage Model}\label{sec:rt_model}
We consider a self-scheduled energy storage in which the storage owner submits bids to sell/purchase energy quantities without specifying the price. The energy storage is assumed to be a price-taker, meaning it does not influence market prices. We formulate the real-time arbitrage problem as follows:
\begin{subequations}\label{ps}
\begin{gather}
\max_{p_t, b_t \in \mathcal{B}} \Big\{\sum_{t\in\mathcal{T}}\hat{\lambda}_t(p_t-b_t) - c p_t\nonumber\Big\} \label{ps:obj}\\
\{p_t,b_t\} \in B ;\forall t\in\mathcal{T} \label{ps:c1}\\
0\leq p_t,b_t \leq P \label{ps:c2}\\
e_t-e_{t-1} = -p_t/\eta + b_t\eta \label{ps:c3}\\
0\leq e_t \leq E \label{ps:c4}\\
p_t=0\;\mathrm{if}\; \lambda_t < 0 \label{ps:c5} 
\end{gather}
\end{subequations}

Here, the set $\mathcal{T}$ includes all intervals within the look-ahead horizon, $\lambda_t$ is the real-time market price, $\hat{\lambda}_t$ is the predicted real-time market price, while $p_t$ is the storage's discharge energy at time $t$, and $b_t$ is the storage's charge energy. The cost parameter $c$ represents discharge costs incurred only during discharge to align with conventional generation cost models. The set $B$ represents a bidding policy, while $\mathcal{B}$ contains all feasible bidding policies. $P$ is the power rating of the storage, $\eta$ is the efficiency, and $e_s$ denotes the state-of-charge (SoC) constraint. 

The objective function maximizes arbitrage profit in the real-time market. The first term reflects revenue from the real-time market, while the second term accounts for physical discharge costs. Equation \eqref{ps:c1} ensures that all storage actions adhere to a feasible bidding policy $B \in \mathcal{B}$, satisfying both market clearing procedures and physical constraints. Equation \eqref{ps:c2} limits the discharge and charge power to the power rating, \eqref{ps:c3} and \eqref{ps:c4} define the SoC dynamics and capacity constraints, where $E$ is the storage capacity. \eqref{ps:c5} prohibits discharging when the price is negative, providing a sufficient condition to prevent simultaneous charging and discharging~\cite{yousuf}.

\subsection{Risk-Averse Real-Time Energy Storage Arbitrage Policy}\label{sec:risk-averse}
We design a risk-averse real-time energy storage arbitrage policy, leveraging the predictions from the previous layers. Specifically, we utilize the point price prediction $\hat{\lambda}$, and a corresponding prediction set $\hat{C}$, defined by the interval $[\hat{\underline{\lambda}}, \hat{\overline{}{\lambda}}]$. The optimization model in \eqref{ps} finds the optimal storage schedule based on point price predictions. However, to incorporate a risk-averse approach, it is important to consider not only the point prediction but also the model's confidence in those predictions.

To achieve this, our risk-averse policy takes into account the prediction sets using the Monte-Carlo Sampling method. Random price samples are drawn from the predicted interval as follows:
\begin{equation}
    \pi_{t}^i = \hat{\underline{\lambda_t}} + (\hat{\overline{\lambda_t}} - \hat{\underline{\lambda_t}}) \cdot U_i, \forall i = 1, 2, \ldots, N.
\end{equation}
where $\pi_{t}^i$ represents a sampled price at time $t$ drawn from the interval ($\hat{C_t}$), and $U_i$ is a uniformly distributed random variable with zero mean and one standard deviation. 

Each sampled instance is input into the real-time arbitrage storage optimization model to obtain optimal decisions. These decisions are aggregated by taking the common decision across samples. This is done by evaluating the charge and discharge decisions for each sample; if all instances recommend charging, the policy chooses to charge; the same applies for discharging, as shown in~\eqref{eq:sd}. Then, the amount of charge or discharge is determined by either a conservative or aggressive setting. In the conservative setting, the minimum value from all samples is selected, while in the aggressive setting, the maximum value is chosen, as detailed in ~\eqref{eq:sca}-\eqref{eq:sda}.
\begin{equation}\label{eq:sd}
    D_t :=
    \begin{cases} 
\text{Charge}_t, & \text{if all} (\hat{b_t^{\pi_{t}^i}}) > 0 \\ 
\text{Discharge}_t, & \text{if all} (\hat{p_t^{\pi_{t}^i}}) > 0 \\ 
0, & \text{o.w.}
    \end{cases}
\end{equation}
\begin{equation}\label{eq:sca}
    b_t :=
    \begin{cases} 
\min(\hat{b_t^{\pi_{t}^i}}), & \text{if}~ D_t = \text{Charge}_t \text{ and  conservative}  \\ 
\max(\hat{b_t^{\pi_{t}^i}}), & \text{if}~ D_t = \text{Charge}_t \text{ and  aggressive}  \\
0, & \text{o.w.}
    \end{cases}
\end{equation}
\begin{equation}\label{eq:sda}
    p_t :=
    \begin{cases} 
\min(\hat{p_t^{\pi_{t}^i}}), & \text{if}~ D_t = \text{Discharge}_t \text{ and  conservative}  \\ 
\max(\hat{p_t^{\pi_{t}^i}}), & \text{if}~ D_t = \text{Discharge}_t \text{ and  aggressive}  \\
0, & \text{o.w.}
    \end{cases}
\end{equation}

\subsection{Full Solution}
Our overall framework is structured as follows: first, we forecast the real-time price using a point prediction model. We then train a conformal control prediction model to construct the prediction sets. Finally, these prediction sets are integrated into our arbitrage model to find the optimal risk-averse storage schedule. Figure~\ref{fig:conformal_diagram} summarizes our framework.

In the first stage, we train a real-time price prediction model to generate point forecasts. We use real-time prices with a look-back window of $L$ as our features. These point forecasts are then passed into the conformal prediction block, which first computes the score function for the calibration data and then trains the conformal prediction model. In the final stage, we perform multiple optimization runs using the prediction sets obtained earlier. The optimal low-risk schedule is determined according to our risk-averse policy. Algorithm~\ref{alg} outlines the full algorithm.

\begin{algorithm2e}
\caption{Conformal Energy Storage Arbitrage}\label{alg}
\KwData{Prepare data according to section 2.1}
\textbf{Require}: prepare the data as described in~\ref{sec:conformal}\\
\textbf{Initialization}:
Set the energy parameters $P, \eta, e_0, E$, set random seed, time horizon $T$, \\
\textbf{Point Prediction}:\\
Train any point prediction model;\\
\textbf{Conformal Prediction}:\\
\For{$i\le M$}{
 \For{$batch \in D_2$}{
 $x \leftarrow \mathbf{X}$\;
 $y \leftarrow \mathbf{\lambda}$\;
 \For{$t \in T$}{
  Compute the signed residuals $(\lambda_t - \hat{\lambda_t})$\;
  Estimate the residuals quantile $q_t$\;
  Predict $\hat{q}_{t+1}$\;
  Estimate $q_{t+1}$ according to \ref{eq:control}\;
  Construct $\hat{C}_t$ according to \ref{eq:prediction_set}\;
 }
 }
 }
 \textbf{Storage Arbitrage}:\\
 Predict the price using trained point prediction model\;
 Construct the prediction sets using conformal control model \eqref{eq:control}\;
 Compute the storage optimal schedule $\hat{p}_t,\hat{b}_t \forall T$ according to \eqref{eq:sd}-\eqref{eq:sda}
\end{algorithm2e}

\section{Case Study}

\subsection{Experiment Setup}

We use the New York City zonal price data, two weeks in 2022 for calibration, and 2023 year for testing. A look-back window of 24 hours and a varied prediction horizon from 1 hour to 24 hours ahead were considered. We consider a 0.5 MW/1 MWh storage, 90\% one-way efficiency, and \$10/MWh fixed discharge cost. For the conformal prediction, we set $\alpha=0.05$, for a desired coverage of $95\%$, reflecting risk-averse settings.
All experiments are run on a high-performance computing cluster with Intel Xeon Platinum 8640Y 2Ghz CPU, 512GB Memory, and 2xNVIDIA A40 GPU. All codes are available on Github\footnote{https://github.com/Alghumayjan/Risk-Averse-Storage-Arbitrage}.

\subsection{RTP Conformal Prediction Results}
We present the results of the real-time electricity price conformal control prediction. To demonstrate the effectiveness of the conformal control, we consider two synthetically generated point-prediction scenarios, defined as:
\begin{equation}
    \hat{\lambda}_{t} = \lambda_{t} + \beta_t 
\end{equation}
where $\beta$ represents a noise error randomly sampled from a Gaussian distribution with zero mean and standard deviation equal to the targeted error. The target error reflects how well the synthetic point predictor is. We consider two point predictors: a \textit{good forecaster} with a targeted error of \$5/MWh, and a \textit{bad forecaster} with a targeted error of \$40/MWh.

\begin{figure}[htbp]
\centerline{\includegraphics[width=0.9\linewidth]{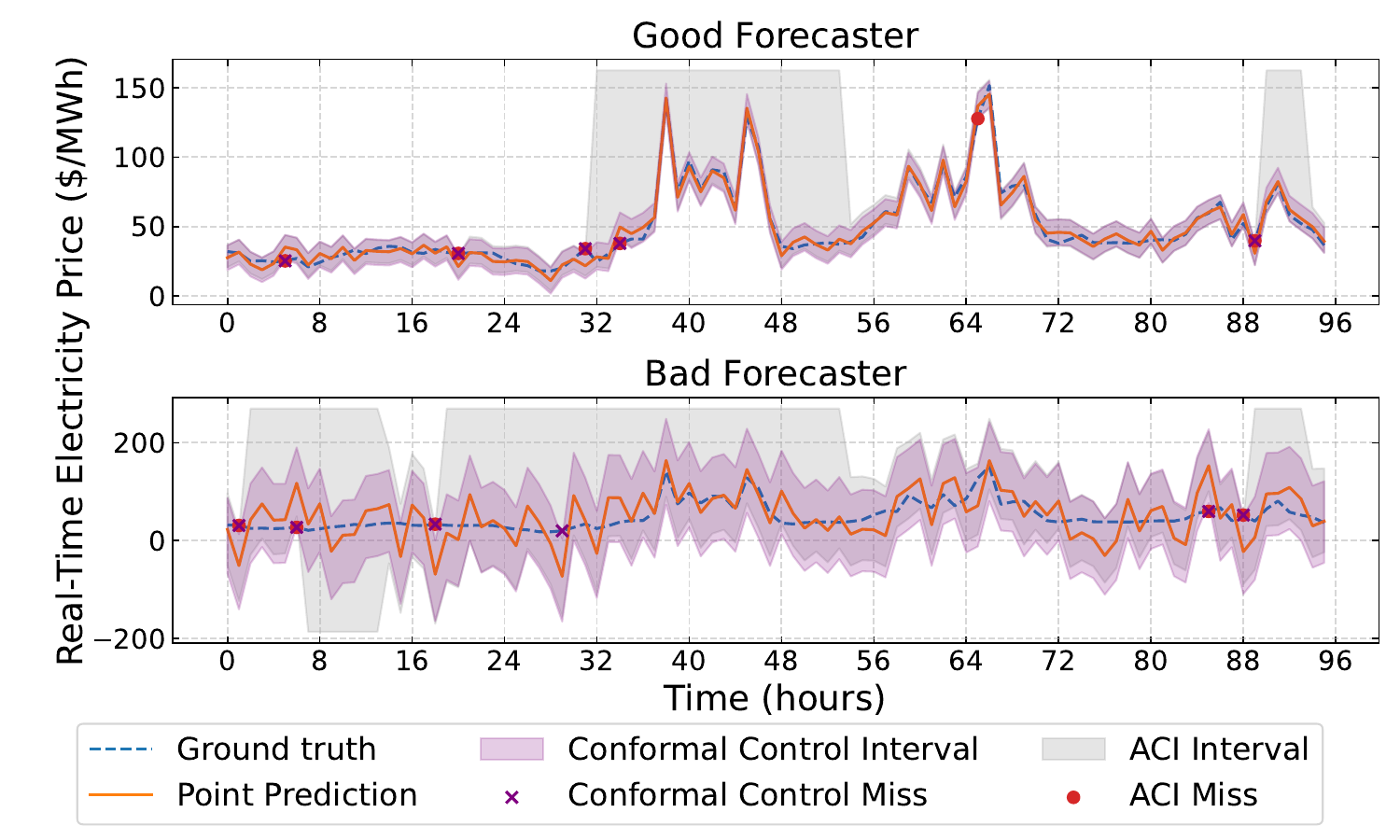}}
\caption{Prediction example for a subset of the testing sample.}
\label{fig:pred_example}
\end{figure}
\begin{figure}[htbp]
\centerline{\includegraphics[width=0.9\linewidth]{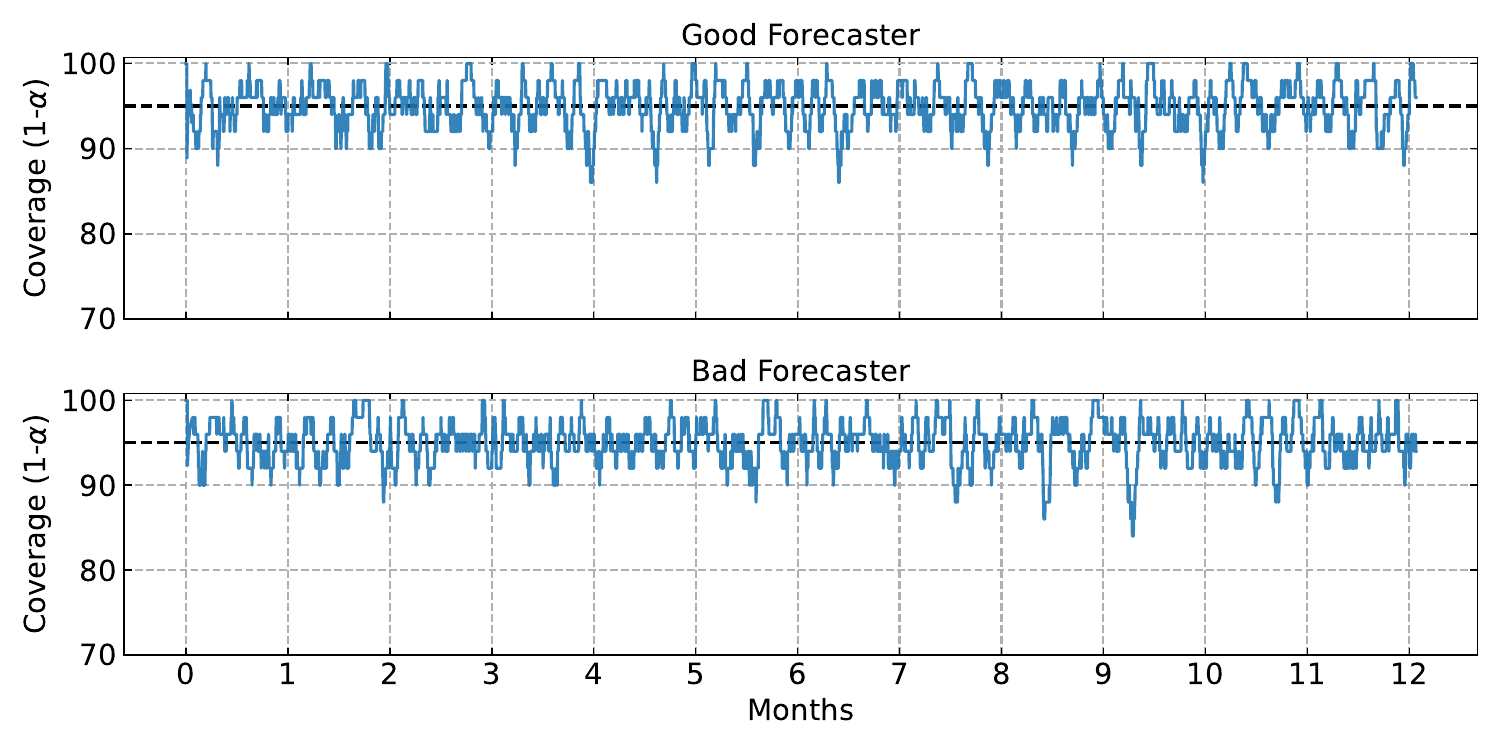}}
\caption{The predicted interval coverage of the ground truth.}
\label{fig:pred_example_cvg}
\end{figure}
Figure~\ref{fig:pred_example} shows a 4-day price example from the testing set, illustrating the predictions made by the conformal control for both scenarios. We compare the conformal control with another conformal prediction model, adaptive conformal inference (ACI) model~\cite{gibbs2021adaptive}. The results show that when we have a good point predictor, the prediction interval remains narrow and consistently covers most of the ground truth. However, if the point predictor is not that good, the interval expands to account for the increased uncertainty. Additionally, we observe that following any missed point, the interval widens and only contracts if it maintains good coverage over a sustained period. This shows how the conformal control corrects itself after each error. In contrast, ACI exhibits similar coverage, but after a series of missed predictions, it starts to predict infinitely wide intervals, as observed in both scenarios. Moreover,~\ref{fig:pred_example_cvg} shows the coverage achieved across the entire testing period, showing that the conformal control maintains the desired coverage of around $95\%$ in both scenarios.

\begin{figure}[htbp]
\centerline{\includegraphics[width=0.9\linewidth]{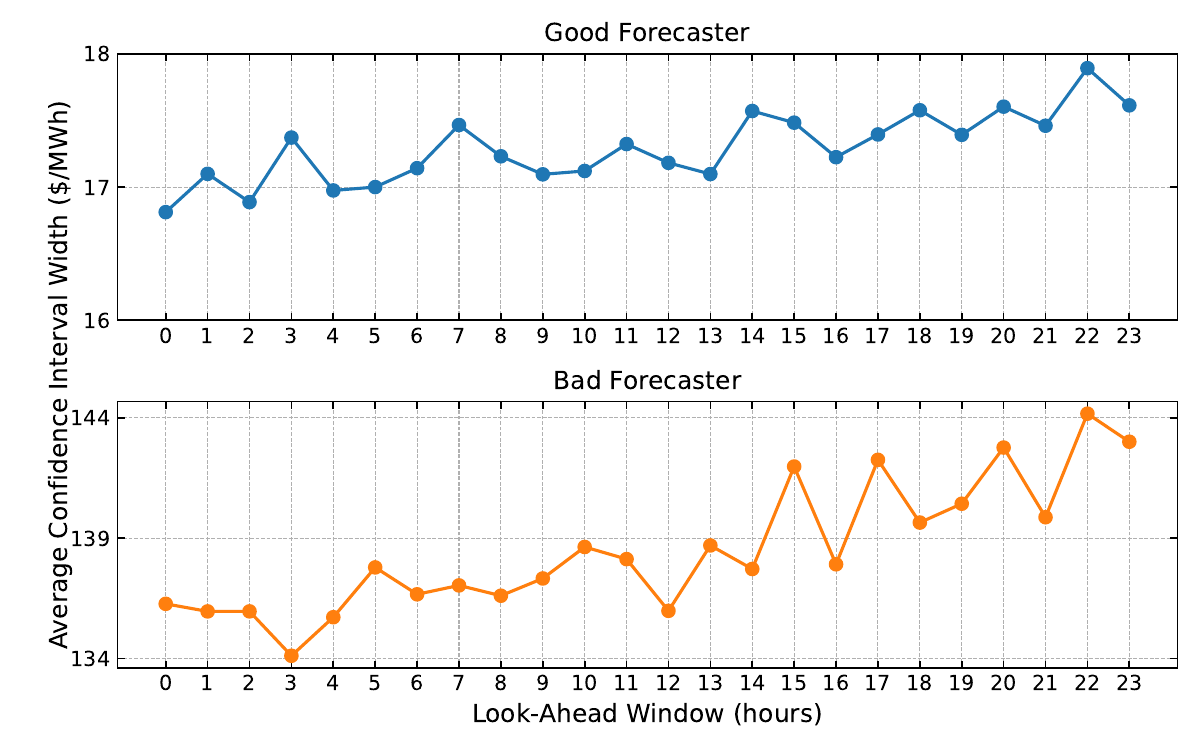}}
\caption{Comparison look-ahead window and the confidence interval.}
\label{fig:horizon}
\end{figure}

In the electricity market, a long forecasting horizon is always desired. The short-term horizon usually lasts from a couple of hours to one day. However, longer horizons mean higher uncertainty. Figure~\ref{fig:horizon} shows how the predicted interval from the conformal control reacts to longer look-ahead windows; as the window gets larger, the average size of the interval will be higher. 

\subsection{Risk-Averse Storage Arbitrage Results}
\begin{figure*}[htbp]
\centerline{\includegraphics[width=0.9\linewidth]{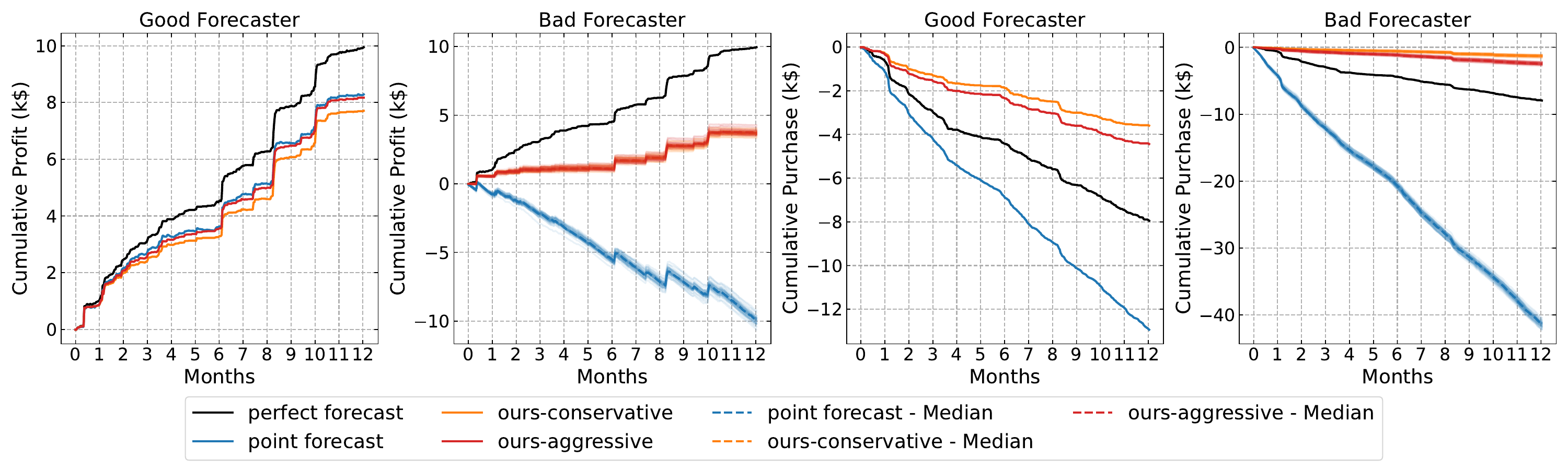}}
\caption{Cumulative daily net profits and purchases following each policy under different scenarios. For the \textit{bad forecaster} scenario, we run multiple simulations with randomly sampled predictions to assess the model's sensitivity.}
\label{fig:profits_purchases}
\end{figure*}
In this subsection, we show the results of our risk-averse policy as outlined in Section \ref{sec:risk-averse}. We consider a prediction horizon of 6 hours. Based on~\ref{fig:horizon} and our experiments, we saw that 6 hours ahead provides a good balance between uncertainty and operational benefits. We run our risk-averse policy with 100 randomized instances within the predicted intervals by our model and use the decisions in two settings: conservative and aggressive.

Figure~\ref{fig:profits_purchases} illustrates the cumulative profits achieved in each scenario. The results show that our policy yields a good profit margin with the \textit{good forecaster} and profiting with the \textit{bad forecaster} where the point forecast case is losing.

The figure also shows the cumulative purchases made by the storage in each scenario. We can see that in all cases, our policy has the lowest amount of purchases, which helped prevent more high losses with the bad predictors. We also note that in the \textit{good forecaster} case, we achieved similar profit as the point forecast case while making less than $35\%$ of its purchases. This helps in avoiding cycling the storage more, which can increase its lifetime. Table~\ref{tab:summary} shows a numerical comparison between the cases. Our risk-averse modes show the lowest purchase amount and negative profit days.

\begin{table}[htbp]
\caption{Results Summary by Forecaster and Case}
\begin{tabular}{ccccc}
\toprule
  &  & \textbf{Neg.} & \textbf{Total} & \textbf{Total} \\
 \textbf{Forecaster}& \textbf{Case} & \textbf{Days} & \textbf{Profit (\$)} & \textbf{Purchases (\$)} \\
\midrule
\multirow{3}{*}{Good} & Perfect Forecast & 54  & 9,947.80  & 7,948.99  \\
\multirow{3}{*}{Forecaster} & Point Forecast   & 148 & \textbf{8,284.20}  & 12,936.43	 \\
 & Ours-Conservative       & 43  & 7,707.87  & \textbf{3,596.60}  \\
 & Ours-Aggressive        & \textbf{35}  & 8,170.12  & 4,438.98  \\
\midrule
\multirow{2}{*}{Bad}  & Perfect Forecast & 54  & 9,947.80  & 7,948.99  \\
\multirow{2}{*}{Forecaster}& Point Forecast   & 320 & -9,902.98 & 41,301.41 \\
\multirow{2}{*}{(Median)}& Ours-Conservative       & \textbf{86}  & 3,697.08	  & \textbf{1,316.20}  \\
& Ours-Aggressive & 104 & \textbf{3,737.29}  & 2,453.90 \\
\bottomrule
\end{tabular}
\label{tab:summary}
\end{table}

\section{Conclusion}
We proposed a risk-averse energy storage price arbitrage utilizing uncertainty quantification for electricity real-time price forecasting via conformal prediction. We built a conformal control model for electricity real-time prices that can work with any underlying point forecast model without assuming any distribution on the dataset. We designed a risk-averse policy for energy storage arbitrage that takes into account the predicted confidence interval and makes the safest decision. Our risk-averse policy has shown to be very robust in avoiding high losses, as it achieved a similar profit margin as the point forecast while making less than $35\%$ purchases.
For the next steps, we plan to improve the conformal control model by including a kernel integrator to better capture the price daily cycles and recurring spike times.

\bibliographystyle{IEEEtran}
\bibliography{references}

\begin{thebibliography}{10}
\providecommand{\url}[1]{#1}
\csname url@samestyle\endcsname
\providecommand{\newblock}{\relax}
\providecommand{\bibinfo}[2]{#2}
\providecommand{\BIBentrySTDinterwordspacing}{\spaceskip=0pt\relax}
\providecommand{\BIBentryALTinterwordstretchfactor}{4}
\providecommand{\BIBentryALTinterwordspacing}{\spaceskip=\fontdimen2\font plus
\BIBentryALTinterwordstretchfactor\fontdimen3\font minus \fontdimen4\font\relax}
\providecommand{\BIBforeignlanguage}[2]{{%
\expandafter\ifx\csname l@#1\endcsname\relax
\typeout{** WARNING: IEEEtran.bst: No hyphenation pattern has been}%
\typeout{** loaded for the language `#1'. Using the pattern for}%
\typeout{** the default language instead.}%
\else
\language=\csname l@#1\endcsname
\fi
#2}}
\providecommand{\BIBdecl}{\relax}
\BIBdecl

\bibitem{wang2017impact}
Y.~Wang, Y.~Dvorkin, R.~Fern{\'a}ndez-Blanco, B.~Xu, and D.~S. Kirschen, ``Impact of local transmission congestion on energy storage arbitrage opportunities,'' in \emph{2017 IEEE Power \& Energy Society General Meeting}.\hskip 1em plus 0.5em minus 0.4em\relax IEEE, 2017, pp. 1--5.

\bibitem{us2023data}
\text{US Energy Information Association}, ``Battery storage in the united states: An update on market trends,'' 2023.

\bibitem{qin2023role}
X.~Qin, B.~Xu, I.~Lestas, Y.~Guo, and H.~Sun, ``The role of electricity market design for energy storage in cost-efficient decarbonization,'' \emph{Joule}, 2023.

\bibitem{lago2018forecasting}
J.~Lago, F.~De~Ridder, and B.~De~Schutter, ``Forecasting spot electricity prices: Deep learning approaches and empirical comparison of traditional algorithms,'' \emph{Applied Energy}, vol. 221, pp. 386--405, 2018.

\bibitem{panapakidis2016day}
I.~P. Panapakidis and A.~S. Dagoumas, ``Day-ahead electricity price forecasting via the application of artificial neural network based models,'' \emph{Applied Energy}, vol. 172, pp. 132--151, 2016.

\bibitem{lehna2022forecasting}
M.~Lehna, F.~Scheller, and H.~Herwartz, ``Forecasting day-ahead electricity prices: A comparison of time series and neural network models taking external regressors into account,'' \emph{Energy Economics}, vol. 106, p. 105742, 2022.

\bibitem{alghumayjan2024energy}
S.~Alghumayjan, J.~Han, N.~Zheng, M.~Yi, and B.~Xu, ``Energy storage arbitrage in two-settlement markets: A transformer-based approach,'' in \emph{Proceedings of the Power System Computation Conference (PSCC) 2024}, Paris, France, 2024.

\bibitem{elmachtoub2022smart}
A.~N. Elmachtoub and P.~Grigas, ``Smart “predict, then optimize”,'' \emph{Management Science}, vol.~68, no.~1, pp. 9--26, 2022.

\bibitem{vovk1999machine}
V.~Vovk, A.~Gammerman, and C.~Saunders, ``Machine-learning applications of algorithmic randomness,'' 1999.

\bibitem{vovk2005algorithmic}
V.~Vovk, A.~Gammerman, and G.~Shafer, \emph{Algorithmic learning in a random world}.\hskip 1em plus 0.5em minus 0.4em\relax Springer, 2005, vol.~29.

\bibitem{angelopoulos2024conformal}
A.~Angelopoulos, E.~Candes, and R.~J. Tibshirani, ``Conformal pid control for time series prediction,'' \emph{Advances in neural information processing systems}, vol.~36, 2024.

\bibitem{renkema2024conformal}
Y.~Renkema, N.~Brinkel, and T.~Alskaif, ``Conformal prediction for stochastic decision-making of pv power in electricity markets,'' \emph{Electric Power Systems Research}, vol. 234, p. 110750, 2024.

\bibitem{yeh2024end}
C.~Yeh, N.~Christianson, A.~Wierman, and Y.~Yue, ``End-to-end conformal calibration for robust grid-scale battery storage optimization,'' in \emph{NeurIPS 2024 Workshop on Tackling Climate Change with Machine Learning}, Vancouver, Canada, 12 2024.

\bibitem{zheng2022arbitraging}
N.~Zheng, J.~Jaworski, and B.~Xu, ``Arbitraging variable efficiency energy storage using analytical stochastic dynamic programming,'' \emph{IEEE Transactions on Power Systems}, vol.~37, no.~6, pp. 4785--4795, 2022.

\bibitem{sioshansi2021energy}
R.~Sioshansi, P.~Denholm, J.~Arteaga, S.~Awara, S.~Bhattacharjee, A.~Botterud, W.~Cole, A.~Cortes, A.~De~Queiroz, J.~DeCarolis \emph{et~al.}, ``Energy-storage modeling: State-of-the-art and future research directions,'' \emph{IEEE transactions on power systems}, vol.~37, no.~2, pp. 860--875, 2021.

\bibitem{yousuf}
Y.~Baker, N.~Zheng, and B.~Xu, ``Transferable energy storage bidder,'' \emph{IEEE Transactions on Power Systems}, 2023.

\bibitem{gibbs2021adaptive}
I.~Gibbs and E.~Candes, ``Adaptive conformal inference under distribution shift,'' \emph{Advances in Neural Information Processing Systems}, vol.~34, pp. 1660--1672, 2021.

\end{thebibliography}

\end{document}